\documentclass[sigplan,screen]{acmart}
\usepackage{algorithm}
\usepackage{algpseudocode}
\usepackage{subfig}
\AtBeginDocument{%
  }

\setcopyright{acmlicensed}
\copyrightyear{2018}
\acmYear{2018}
\acmDOI{XXXXXXX.XXXXXXX}

\begin{document}

\title{Analysis of Different Algorithmic Design Techniques for Seam Carving}

\author{S Muhammad Ali}
\email{sn07590@st.habib.edu.pk}
\affiliation{%
  \institution{Habib University}
  \city{Karachi}
  \country{Pakistan}
}

\author{Owais Aijaz}
\email{oa07610@st.habib.edu.pk}
\affiliation{%
  \institution{Habib University}
  \city{Karachi}
  \country{Pakistan}
}

\author{Yousuf Uyghur}
\email{mu07486@st.habib.edu.pk}
\affiliation{%
  \institution{Habib University}
  \city{Karachi}
  \country{Pakistan}
}


\begin{abstract}
 Seam carving, a content-aware image resizing technique, has garnered significant attention for its ability to resize images while preserving important content. In this paper, we conduct a comprehensive analysis of four algorithmic design techniques for seam carving: brute-force, greedy, dynamic programming, and GPU-based parallel algorithms. We begin by presenting a theoretical overview of each technique, discussing their underlying principles and computational complexities. Subsequently, we delve into empirical evaluations, comparing the performance of these algorithms in terms of runtime efficiency. Our experimental results provide insights into the theoretical complexities of the design techniques. 
\end{abstract}



\begin{teaserfigure}
  \includegraphics[width=\textwidth]{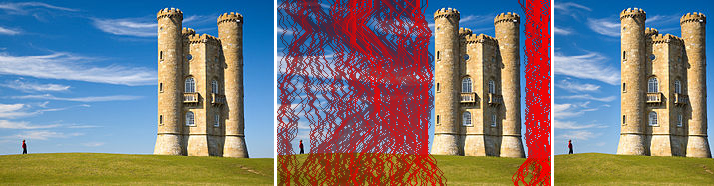}
  \caption{Left to right: Original Image to be resized, List of possible seams, Final Image}
  \label{fig:teaser}
\end{teaserfigure}

\maketitle

\section{Introduction}
Resizing images is a fundamental operation in image processing, facilitating a myriad of applications across various domains. Traditionally, methods like cropping and scaling have been employed for this purpose. However, these techniques often lead to undesirable outcomes such as loss of data or visible distortions in image content. As digital imagery becomes increasingly prevalent, the demand for content-aware resizing methods has surged, driving researchers to explore innovative solutions.

A notable breakthrough in this realm is the seam carving algorithm introduced by Avidan and Shamir \cite{Avidan}. Seam carving operates by selectively removing or adding pixel paths, known as seams, to resize images while preserving important content and foreground objects. At its core, seam carving identifies optimal paths of pixels—seams—based on an image energy function. By iteratively removing the least significant pixel paths, or adding new paths between less important ones, seam carving mitigates the shortcomings of traditional resizing methods, by selectively preserving crucial image elements while reducing or expanding its dimensions, offering superior visual quality and content retention. However, the computational complexity of seam carving, especially in establishing energy maps and identifying optimal seams, can impact its operational speed compared to traditional resizing methods like clipping or scaling.

This report explores various aspects of seam carving, including its foundational principles, design techniques, and real-world applications. Additionally, it examines performance enhancements and alternative approaches proposed in the literature to address challenges such as computational complexity and preservation of image structure. We will talk about four distinct algorithmic paradigms: dynamic programming, greedy, brute force, and dynamic GPU. Each approach brings its own set of design techniques and computational considerations, all with the shared goal of optimizing the seam carving process. Additionally, we'll conduct empirical analyses to compare algorithm performance across different input sizes, shedding light on their scalability and real-world efficiency.

\section{Background}

Image resizing is a common feature in many image processing applications, typically achieved by uniformly adjusting the image to a desired size. However, there's a growing interest in image retargeting, a process that aims to alter the image dimensions while preserving important features. These features can be identified either top-down or bottom-up. Top-down methods utilize tools like face detectors to pinpoint crucial regions in the image, while bottom-up approaches rely on visual saliency methods to create a saliency map highlighting significant areas \cite{voila}. Once this map is generated, cropping can be applied to focus on the most crucial parts of the image. Various studies have explored different aspects of image retargeting, including automatic thumbnail creation, adapting images for mobile devices, and intelligent cropping using techniques like eye tracking and composition rules \cite{suh}.

Traditionally, image resizing and cropping operations have been the primary means of altering image sizes. However, researchers have proposed more advanced techniques, such as non-linear, data-dependent scaling. For instance, Liu and Gleicher introduced a method for image and video retargeting that identifies Regions of Interest (ROI) and applies a novel Fisheye-View warp to preserve these areas while adjusting the rest of the image or video \cite{liu}.

Another approach involves decomposing images into background layers and foreground objects, allowing for the removal and resizing of non-essential elements before reinserting important regions \cite{suh}. Additionally, Gal et al. introduced a solution for warping images into arbitrary shapes while maintaining user-specified features through Laplacian editing techniques..

Seam carving, a technique used in various image editing applications, involves finding optimal seams to combine parts of images seamlessly. It has been applied in systems like Digital Photomontage and AutoCollage \cite{Avidan}. However, most of these methods focus on image editing rather than retargeting.

Various algorithms and techniques are used to compute seams, including Dijkstra's shortest path algorithm, dynamic programming, and graph cuts. Seam carving has also found applications in texture synthesis and object removal, where it helps in generating larger texture images and filling in missing parts of images respectively \cite{voila}.

Overall, the field of image manipulation and retargeting has seen significant advancements through the development of innovative algorithms and methodologies aimed at preserving important features while altering image dimensions to meet specific requirements.

\subsection{Enhancements and Applications}

Seam carving, originally introduced as an innovative technique for content-aware image resizing, has undergone significant enhancements and found diverse applications across various domains in image processing and computer vision. In this section, we explore some of the key enhancements made to the seam carving algorithm and its wide-ranging applications.

\subsubsection{Forward Energy Calculation}
Forward energy is a modification to the traditional seam carving algorithm that considers the energy created from the new neighbors after a seam is removed. In the original algorithm, the accumulated cost matrix is constructed based solely on the existing energy values. However, with forward energy, the algorithm incorporates the additional energy introduced by the new neighbors created after seam removal. This forward energy is added during the computation of the accumulated cost matrix, influencing the selection of optimal seams.

The benefits of forward energy usage are evident in preserving the structural integrity and visual coherence of the resized image. By accounting for the energy contribution of new neighbors, forward energy helps maintain important visual features and prevent distortions, as demonstrated by examples showcasing the preservation of angular shapes and overall image quality.

\subsubsection{Non-linear Seam Carving}
Traditional seam carving methods employ linear resizing, which may not always provide optimal results, especially in complex images with non-uniform content distribution. Non-linear seam carving techniques address this limitation by allowing for more flexible resizing operations, adapting to the image content in a nonlinear fashion.

\subsubsection{Image Enlargement}
Seam insertion complements seam removal by allowing the insertion of new seams into an image to either enlarge it or restore content that was previously removed. This process involves performing seam removal in reverse to identify the coordinates and order of minimum seams that were originally removed. The original image remains unaltered during this phase, as the purpose is solely to record the coordinates and order of the minimum seams. Subsequently, the identified coordinates are used to insert new seams into the original image in the same order they were removed. The pixel values of the new seams are typically derived from an average of the top/bottom or left/right neighbors.

\subsubsection{Image Editing}
Seam carving techniques have found applications beyond resizing. It can be integrated into interactive image editing tools, allowing users to manipulate images while receiving real-time feedback on the preservation of important features. Seam carving can assist in retouching and restoring damaged or degraded images. By selectively removing or adding seams, imperfections can be corrected or obscured, resulting in visually improved images \cite{liu}. 

\subsubsection{Texture Synthesis}
Seam carving plays a role in texture synthesis by intelligently stitching together texture patches to create larger, seamless textures. This is beneficial in graphics and design applications where high-resolution textures are required.

\subsubsection{Image Composition}
In applications like digital photomontage and collage creation, seam carving aids in seamlessly blending multiple images together. By finding optimal seams to combine image fragments, seam carving contributes to creating visually appealing compositions \cite{suh}.

\subsubsection{Image Retargeting}
Seam carving facilitates the adaptation of images to different aspect ratios or screen sizes without losing critical information. This is valuable for optimizing images for display on various devices, such as smartphones, tablets, and computer screens.

\subsubsection{Object Removal}
Seam carving can be employed to remove unwanted objects or elements from images seamlessly. By identifying and removing seams containing undesirable content, seam carving effectively cleans up images without leaving noticeable artifacts. This process begins with the creation of a binary mask during preprocessing, outlining the area to be removed. The algorithm then evaluates whether it is more efficient to remove the object by removing top to bottom seams or left to right seams, based on the masked region's area. When generating the image's energy map, areas under the masked region are assigned a very high negative value, ensuring that the minimum seam will pass through the masked region. Seam removal is performed as usual, with the additional step of updating the mask to reflect the removal of the minimum seam. Once the masked region has been completely removed, seam insertion is employed to restore the image to its original dimensions

\section{Preliminary}
Content-aware image resizing has emerged as a powerful technique for adjusting the dimensions of images while preserving their essential features and visual coherence. Now to achieve this, the approach is to remove pixels in a judicious manner. Therefore, the question is how to chose the pixels to be removed? The following subsections talk about the how the energy mapping and the seam-carving algorithm as whole works. This flowchart from \cite{IJASCE} shows the algorithm in practice.

\begin{figure}[hbtp]
  \centering
  \includegraphics[width=\linewidth]{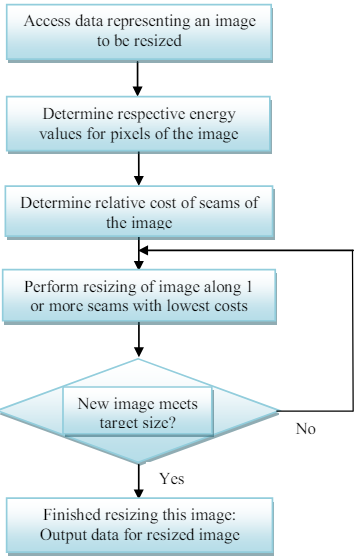}
  \caption{Flow chart to reduce image size}
\end{figure}

Intuitively, our goal is to remove unnoticeable pixels that blend with their surroundings. This leads to the simple energy function shown below. However, there are other energy functions as well which would be described in the following subsections. Here's how the energy mapping works:

\subsection{Energy Mapping}
The concept of energy in an image refers to the importance of each pixel in conveying the overall visual content. High-energy areas typically contain significant visual features like edges, textures, or high contrast regions, while low-energy areas consist of relatively uniform or less important parts of the image. Energy is computed based on gradients within the image. Gradients represent the rate of change of pixel intensity in various directions. Typically, the energy at each pixel is calculated by summing the magnitudes of gradients in both the horizontal and vertical directions. Some common types of energy functions defined in the paper \cite{wedler} are:

\begin{itemize}
\item {\texttt{e1}}: It is the sum of the magnitude of x and y partial derivatives. This is the simplest energy function and gives good results as well. 
\item {\texttt{HoG}}: It is \texttt{e1} divided by a histogram of gradients measure. The intention is to attract a seam to an edge in the picture but to not cross it. This is done by dividing by the max of a histogram of gradients measure to reduce the cost near an edge. But then the cost is multiplied by the e1 energy, so it is high at an edge.
\item {\texttt{e2}}: It uses the second derivative combined with e1.
\item {\texttt{entropy}}: It calculates the entropy along a seam. This is a very different approach because it is not based on gradients, but on the distribution of color and brightness. The optimal seam is one with the greatest uniformity, so entropy is minimized.
\end{itemize}

Once the energy values for all pixels are computed, they are used to create an energy map, which is essentially a matrix representing the energy values at each pixel location in the image. The energy map provides a visual representation of which areas in the image contain the most important visual information. Figure 3 shows a sample energy map,

\begin{figure}[hbtp]
  \centering
  \includegraphics[width=\linewidth]{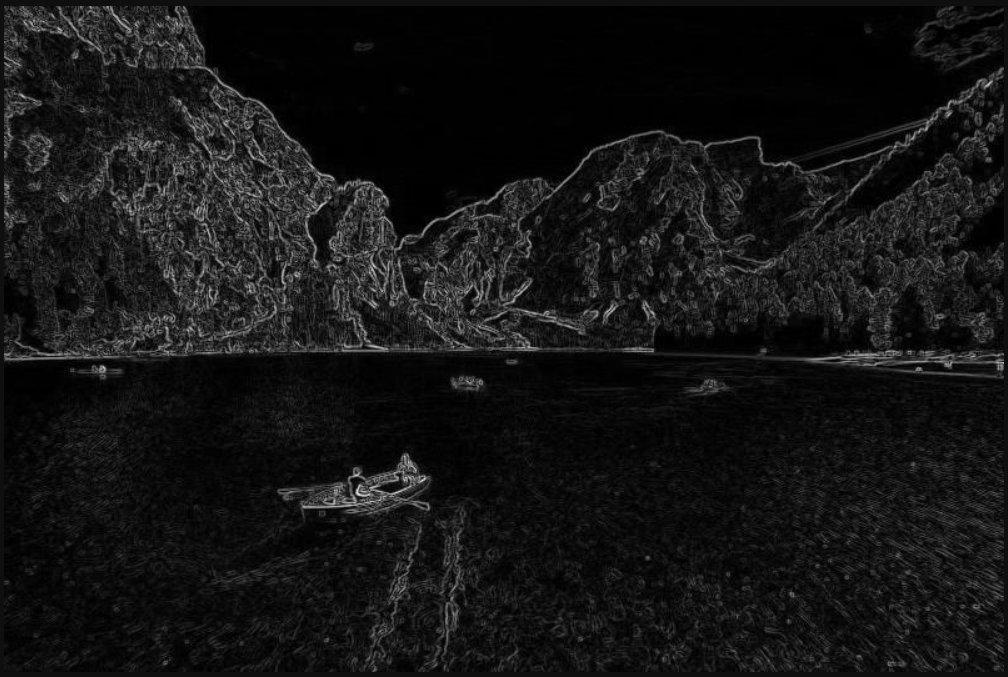}
  \caption{Visualization of an energy map}
\end{figure}

\subsection{Seam Identification}
After the energy map is generated, the seam carving algorithm identifies seams. A seam in the context of image resizing refers to a connected path of pixels traversing the image from one side to another. Seams can be either vertical or horizontal, depending on the desired resizing direction. Seams are a connected path which means that they can only deviate by $\pm$ 1 px between a row or a column. These seams are chosen based on their low cumulative energy, indicating that they pass through areas of relatively less importance in the image. This is the part where different design techniques come in. To efficiently find the optimal seams with the lowest cumulative energy, dynamic programming techniques are often employed. Dynamic programming allows the algorithm to efficiently compute and store intermediate results, enabling it to find the optimal seam without redundantly recomputing energy values. However, we have implemented the same algorithm with different design techniques and contrasted them in the later sections.

\begin{figure}[hbtp]
  \centering
  \includegraphics[width=\linewidth]{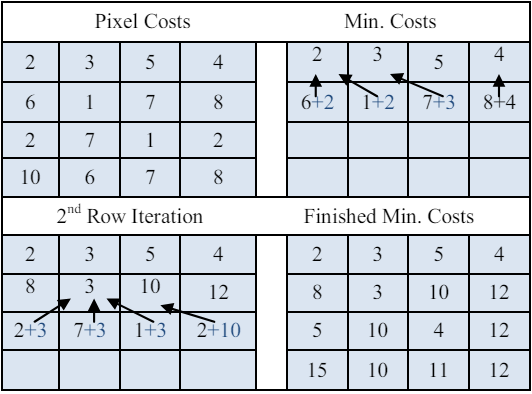}
  \caption{Forming minimum cost table \cite{IJASCE}}
\end{figure}

Figure shows an example for finding an optimal seam.
A minimum cost table in which energy values of each pixel
is calculated to find minimum-cost path. To form this table,
it is assumed that we are intended to find seams from up to
down, top three neighbor pixels of the current pixel are
benefited. The least pixel value of those three neighbors is
added to current pixel’s value and the cost of the
corresponding pixel is obtained as that sum. Operation goes
as follows:
Minimum seam is determined after the minimum cost
table is established. The pixel which has the minimum value
at the bottom row of the table will be the bottom pixel of the
optimal seam. Then, a pixel which has the minimum value
among above three neighbors of that pixel is chosen as a new
part of the seam and the pixel selection goes so. 

\begin{figure}[hbtp]
  \centering
  \includegraphics[width=\linewidth]{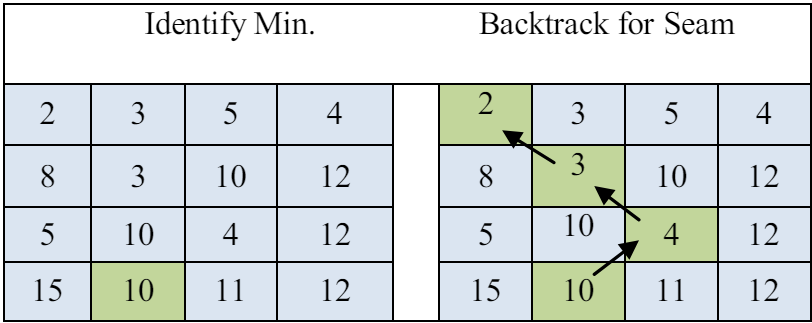}
  \caption{Determining minimum (optimal) seam \cite{IJASCE}}
\end{figure}

\subsection{Removal of Seams}
Once the optimal seam is identified, the seam carving algorithm proceeds to remove or add the seam, depending on the desired resizing operation. When reducing the image size, seams with the lowest energy are removed, effectively eliminating less visually significant pixels. Conversely, when enlarging the image, new seams are added to introduce additional pixels, ensuring that the resized image maintains its essential features and visual coherence.

\hfill \break All the pixels in each row after the pixel to be removed are shifted over one column (or row). Finally, the width (or height) of the image has been reduced by exactly one pixel. By iteratively identifying and removing low-energy seams from the image, seam carving achieves content-aware resizing that maintains the integrity of important visual elements, resulting in more visually pleasing results compared to traditional resizing techniques.

\section{Design Techniques}
This section covers multiple approaches to the seam carving problem. The Seam carving algorithm has been implemented and tested using four different design techniques which include the naive brute force approach, greedy approach, dynamic programming approach and a gpu accelerated approach. This section first makes some assumptions and then describes the approaches in details and their theoretical complexity analysis.

For simplicity we will talk about vertical seams (column-wise) throughout our analysis, since the problem of finding horizontal seams can be reduced to finding minimum vertical seam(s) without loss of generality. We can either replace rows with columns and vice versa in the algorithm or just rotate our image by 90 degrees to find horizontal seams.

\subsubsection{Assumptions}
Let $I$ be an image of dimensions $n$x $m$ where $n,m$ is the number of rows and columns respectively. Let $m'$ be the desired number of columns that we want our image to be reduced to, so the total number of times seam carving needs to run is $c' = (1-s)*m = s'm$ times where $s$ is the ratio of no of columns of resulting image to the original image. Let $E$ be the $m$x $n$ array where every entry $e(i,j)$ corresponds to energy of every pixel calculated using any energy function.

\subsection{BruteForce Seam Carving}
The most naive approach to solving the problem of finding the minimum vertical seam in an image is to iterate over all the seams column wise and calculate the total energy for every possible seams which is possible from top to bottom. The recurrence is as under
\[ M(i,j) = e(i,j) + min(M(i-1,j-1), M(i-1, j), M(i-1, j+1) \]

\renewcommand{\algorithmicrequire}{\textbf{Input:}}
\renewcommand{\algorithmicensure}{\textbf{Output:}}
\renewcommand{\algorithmicprocedure}{\textbf{Function}}

\begin{algorithm}
\caption{Seam carving Brute force algorithm}\label{alg:cap}
\begin{algorithmic}
\Require $E$, $n$, $m$
\Ensure $S$ (minimum Seam)
\State $minCost \gets inf$
\State $S \gets []$

\Procedure{SeamHelper}{$row, col, cost=0, Seam=[]$}
        \If{row = 0}
            \If{cost < minCost}
                \State $minCost \gets cost$
                \State $S \gets Seam$
                \State \Return
            \EndIf
        \Else{.}
            \For{$b \gets col-1 \text{ to } col+1$}
                \State \Call{SeamHelper}{$row-1,b,sum+E[row,col],Seam+[col]$}
            \EndFor
            \State \Return
        \EndIf
\EndProcedure

\For{$col \gets 1 \text{ to } m$}
    \State \Call{SeamHelper}{$r, col$}
    
\EndFor

\State \Return $S$

\end{algorithmic}
\end{algorithm}

\subsubsection{Time Complexity}
Let $T(n,m)$ be the time to determine the minimum seam for an image $I$ for size $nxm$ as described above. Then from the algorithm above, we see that function SeamHelper is called $m$ times. Hence,
\[ T(n,m) = \sum_{c=1}^m T(n,c)\]

Where $T(n,c)$ is the time to find the minimum seam starting from the column $c$. For the $cth$ column, we make three recursive calls
\[ T(n,c) = 3T(n-1, c)\]
Let $T(n,c) = r^n$
\[r^n = 3r^{n-1}\]
\[r^{n-1}(r-3) = 0\]
\[ r = 3\]
\[ T(n,c) = 3^n \forall c\]
Now our equation becomes
\[  T(n,m) = \sum_{c=1}^m 3^n\]
\[  T(n,m) = m3^n\]
\[ = \Theta(m3^n)\]

If $m=n$ then complexity is $\Theta(n3^n)$

Hence the time complexity is polynomially exponential. Its worst and best case are also $O(m3^n)$ and $\Omega(m3^n)$ since algorithm has to iterate over all rows and columns in every case. This is the time complexity for remove one seam. Let $C(n,m,s')$ be the cost of removing  removing $s'm$ seams for an $nxm$ image then, 
\[ C(n,m,s') = \sum_{k=0}^{s'm-1} T(n,m-k)\]
\[ = m3^n + (m-1)3^n + ... + (m-s'm-1)3^n\]
\[ = m(m3^n) - 3^n(1 + 2 + .. + s'm-1) \]
\[ = m^23^n - 3^n(s'm(s'm-1)/2) \]
\[ = m3^n (m - s'(s'm-1)/2)\]
If $s'$ is a constant then the overall cost to scale an image to $s$ is,
\[  = \Theta(m3^n * O(m)) \]

This shows that the brute force approach is very time expensive and we need a better design technique to solve this problem.

In general, the overall time for a scale $s$ for all algorithms will be,
\[ = \Theta(T(n) * O(m))\]
In the best case, $s'm = 0$ we get $\Omega(T(n,m))$ \\
In the worst case, $s'm = m$ and we get $O(mT(n,m))$ \\

\subsection{Greedy Seam Carving}
The greedy seam carving approach, at each row, selects the pixel with the lowest cost $e(i,j)$ given that it is a neighbouring pixel of the pixel in the previous row. For the first iteration, it selects the minimum pixel from the $nth$ row. Then it successively selects the minimum out of the 3 neighbouring pixels in the row above it and so on. Let $Sc(i)$ be the cost of minimum seam for an image $I$ at $ith$ row then, 
\[ Sc(i) = e(i,j) + min(E[i-1,j], E(i-1, j+1), E[i-1, j-1)) \]
Here, we can observe that once a column $j$ is selected at the beginning of the algorithm, the algorithm just iterates over $n$ rows to pick the local minimum.

However, since selecting local minimum does not guarantee a correct solution in every case, this algorithm does not produce quality results as compared t =o other approaches

\begin{algorithm}
\caption{Seam carving Greedy algorithm}\label{alg:cap}
\begin{algorithmic}
\Require $E$, $n$, $m$
\Ensure $S$ (minimum Seam)

\State $start \gets \Call{argmin}{E[n, :]}$
\State $cost \gets E[n, start] $
\State $S \gets [start]$
    
    \For{$i \gets n-1$ \textbf{to} $1$}
    \State $idx \gets \Call{argmin}{E[i, start-1], E[i, start], E[i, start+1]}$
            \State $start \gets start + idx - 1$
            \State $S \gets S + [start]$
            \State $cost \gets cost + E[i, start]$
    \EndFor

\State \Return $S$
\end{algorithmic}
\end{algorithm}

\subsubsection{Time Complexity}
The Greedy algorithm selects a minimum from $nth$ row in first stage and then select the minimum pixel out of three adjacent columns for every subsequent row. Hence, it only solves one subproblem at every step.
Let $T(n,m)$ be the time taken to find a seam for an image $I$ of n rows and m columns then
$T(n,m) = m + R(n)$
$m$ is the time taken to find the pixel with minimum energy in the last ($nth$) column and $R(n)$ is the time to iterate over n rows to find minimum seam. It defined as
\[ R(n) = R(n-1) + c \]
Here $c$ is the constant time taken to find the minimum of the energy of three adjacent pixels.
By back substitution,
\[ R(n) = (R(n-2) + c) + c\]
\[ R(n) = R(n-3) +  3c\]
\[ R(n) = R(0) + nc\]
\[R(n) = 0 + nc\]
\[ R(n) = \Theta(n)\]

Now substituting in our equation, we get
\[ T(n,m) = \Theta(m) + \Theta(n) \]
\[ T(n,m) = \Theta(m + n)\]
\[T(n,m) = O(max(m,n))\]

If $n=m$ then the cost becomes $\Theta(n)$. The best and worst case complexities are also the same. \\
Hence the greedy solution give a minimum seam in linear time. However, it does not guarantee an optimal solution. 

 For $s'm = s'n$ seams, we can use the cost function $C$ given in the previous part. If $s$ is a constant, the overall time complexity of the whole procedure we get will be in the order of $O(nT(n)) = O(n^2)$.

\subsection{DP Seam Carving}
The original solution proposed by \cite{Avidan} for the seam carving algorithm uses dynamic programming since it gives an optimal solution in a reasonable time. The recurrence relation is given as
\[ M(i,j) = e(i,j) + min(M(i-1,j-1), M(i-1, j), M(i-1, j+1) \]
\subsubsection{Overlapping Subproblems}
The seam carving problem satisfies the overlapping subproblems property which can be observed in its recurrence equation 
\subsubsection{Optimal Substructure}
The Seam carving problem also follows the optimal substructure property since the an minimum cost of a sub-seam (i.e a part of the seam) gives the minimum cost of that seam. Suppose $M(n,j)$ is the optimal cost of the seam ending at $jth$ column. Let $C$ be the cost of minimum sub-seam upto $n-1$ rows adjacent to $M(n,j)$ then, 
\[M(n,j) = e(m,j) + C\]
Suppose now $C'$ is the cost of another sub-seam adjacent to $M(n,j)$ such that $C' < C$, so
\[e(m,j) + C' < M(n,j) = e(m,j) + C\]
Hence it follows the optimal substructure property.

Therefore a dynamic solution will exist for the problem


\begin{algorithm}
\caption{Seam Carving Dynamic algorithm \cite{kk}}\label{alg:cap}
\begin{algorithmic}
\Require $E$, $n$, $m$
\Ensure $S$ (minimum Seam)

\State $S \gets []$
\State $M \gets \text{initialize a 2d tabular table}$
\State $B \gets \text{initialize a 2d table for backtracking indices}$

 \For{$i \gets 2$ \textbf{to} $n$}
        \For{$j \gets 1$ \textbf{to} $m$}
            \State $idx \gets \Call{argmin}{M[i - 1, j - 1],M[i-1, j],M[i-1, j+1]}$
            \State $B[i, j] \gets idx + j-1$
            \State $M[i,j] \gets E(i,j) + M[i - 1, idx + j-1]$
        \EndFor
    \EndFor

\State $j \gets \Call{argmin}{M(n, :)}$
\For{$k \gets n \text{ to } 1$}
    \State $S \gets S + B[k,j]$
    \State $j \gets B[k,j]$
\EndFor
    
\State \Return $S$

\end{algorithmic}
\end{algorithm}

\subsubsection{Time Complexity}
The algorithm above runs a loop for $n$ times. For every $ith$ iteration from $2$ to $n$ it runs $m$ times. In every $(i,j)th$ iteration, the algorithm performs a constant work. When the loops finishes, the algorithm finds the minimum cost out of $m$ such costs. It then runs $n$ times to stores indices of the seam. Let $T(n,m)$ be the total cost then,
\[ T(n,m) = n * m *c + m + n\]
\[ T(n,m) = O(nm + m + n)\]
\[ T(n,m) = O(nm)\]

if $n=m$ then the cost becomes $O(n^2)$. Hence the dynamic approach takes quadratic time to run the algorithm. \\
Doing it $s'm = s'n$ times will give us an overall cost of $O(nT(n)$ = $O(n^3)$

\subsection{Parallel Seam Carving}
Leveraging the Parallel computing capabilties of Graphic Processor Units, we can parallelize our dynamic solution to make it more efficient. We observe in the dynamic programming approach that while iterating over the rows, at some $ith$ row, the computations for $m$ columns are independent of the others. Hence they can be performed in parallel. The GPU approach is as follows

\begin{algorithm}
\caption{Seam carving GPU algorithm}\label{alg:cap}
\begin{algorithmic}
\Require $E$, $n$, $m$
\Ensure $S$ (minimum Seam)

\State $S \gets []$
\State $M \gets \text{initialize a 2d tabular table}$
\State $B \gets \text{initialize a 2d table for backtracking indices}$

 \For{$i \gets 2$ \textbf{to} $n$}
        \For{$j \gets 1$ \textbf{to} $m$ \textbf{in parallel}}
            \State $idx \gets \Call{argmin}{M[i - 1, j - 1],M[i-1, j],M[i-1, j+1]}$
            \State $B[i, j] \gets idx + j-1$
            \State $M[i,j] \gets E(i,j) + M[i - 1, idx + j-1]$
        \EndFor
    \EndFor

\State $j \gets \Call{argmin}{M(n, :)}$
\For{$k \gets n to 1$}
    \State $S \gets S + B[k,j]$
    \State $j \gets B[k,j]$
\EndFor
    
\State \Return $S$

\end{algorithmic}
\end{algorithm}

\subsubsection{Time Complexity}
Since the $j$ loop is now parallelized, it takes an overall constant time to perform computation for $j$ loop in parallel. Hence now the time complexity becomes

\[ T(n,m) = n*c + n + m\]
\[ T(n,m) = O(n)\]

Hence using GPU, the time complexity for the dynamic solution is reduced to linear time. Finding seams $s'm = s'n$ times will incur a cost of $O(n^2)$ total time.

The results in practice however, depend on the implementation of the GPU algorithm since there hidden costs involved in GPU algorithms which include transferring data from GPU to CPU memory and vice versa, invoking GPU etc. An efficient GPU implementation minimizes such hidden costs hence, bringing the empirical time complexity closer to theoretically linear complexity.

\section{Hardware Implementation and Experiments}

We perform our experiments on images of which have small difference in the ratio of width to height. i.e they are closer to $n$ x $n$ square images. 

\subsubsection{Hardware and implementation specifications}

We implement the above given algorithms in python environemnt. We analyze the empirical runtime of our code by running it on Google Colab. It uses the following hardware. \\
GPU : \textbf{Nvidia T4} \\
CPU:   \textbf{Intel Xeon CPU 2 vCPUs (virtual CPUs)} \\
Memory:   \textbf{13GB} \\
Moreover, it uses python version 3.10.12 on Ubuntu 22.04.3 LTS

\subsection{Runtime For different sized images}
This section presents the resulting plots for an image of different sizes. We use google Gemini AI to generate an $n \times n$ image in figure 6. We then use different sized version of the same image to evaluate the runtime of our algorithms. Specifically, we use this image at size $180 \times 180$, $360 \times 360$, $480 \times 480$. $720 \times 720$, $1080 \times 1080$ for all algorithms except bruteforce since this approaches takes polynomially exponential time to run. We first present the time for whole procedure of carving an image at scale $s = 0.5$ as well as the time for running only the seam carving procedure once. We show that the overall time complexity matches that of the cost $C(n,m,s)$ and the time for running the seam carving algorithm matches the theoretical time $T(n,m)$ calculated above for each design technique.

\begin{figure}
    \centering
    \includegraphics[width=0.75\linewidth]{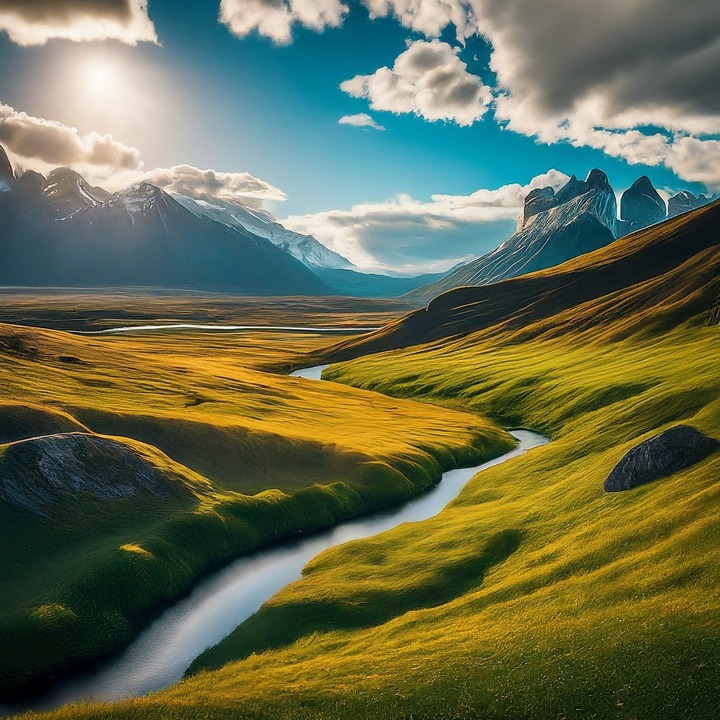}
    \caption{AI generated image used for runtime analysis}
    \label{fig:enter-label}
\end{figure}

\subsection{Brute Force Algorithm}
Figure 7, and 8 shows the runtime for bruteforce algorithm for images of dimensions $2 \times 2$, $5 \times 5$, $7 \times 7$, $10 \times 10$, $12 \times 12$. Going above these dimensions is not feasible to run due to its complexity. The results of the plots are the same as our theoretical analyses i.e $C(n,n,s) = O(n^23^n)$ and $T(n,n) = O(n3^n)$

\begin{figure}[hbtp]
    \centering
    \includegraphics[width=0.75\linewidth]{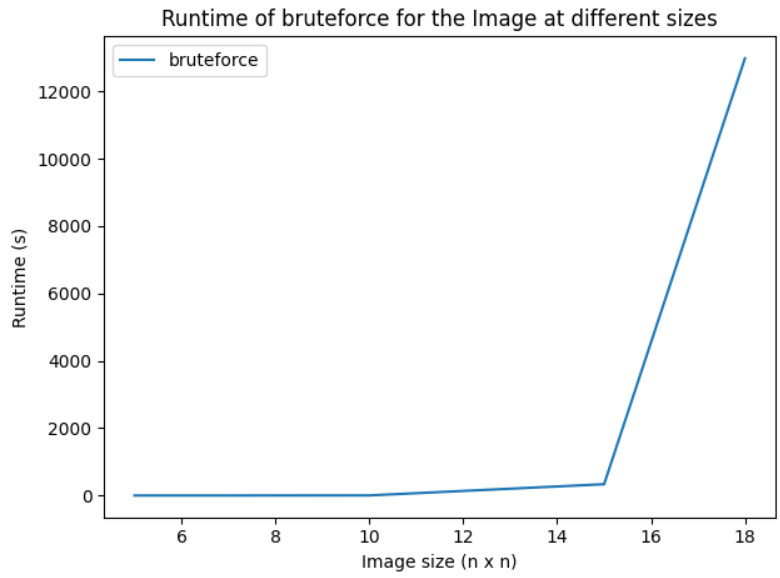}
    \caption{Overall time complexity of the brute force algorithm}
    \label{fig:enter-label}
\end{figure}

\begin{figure}[hbtp]
    \centering
    \includegraphics[width=0.75\linewidth]{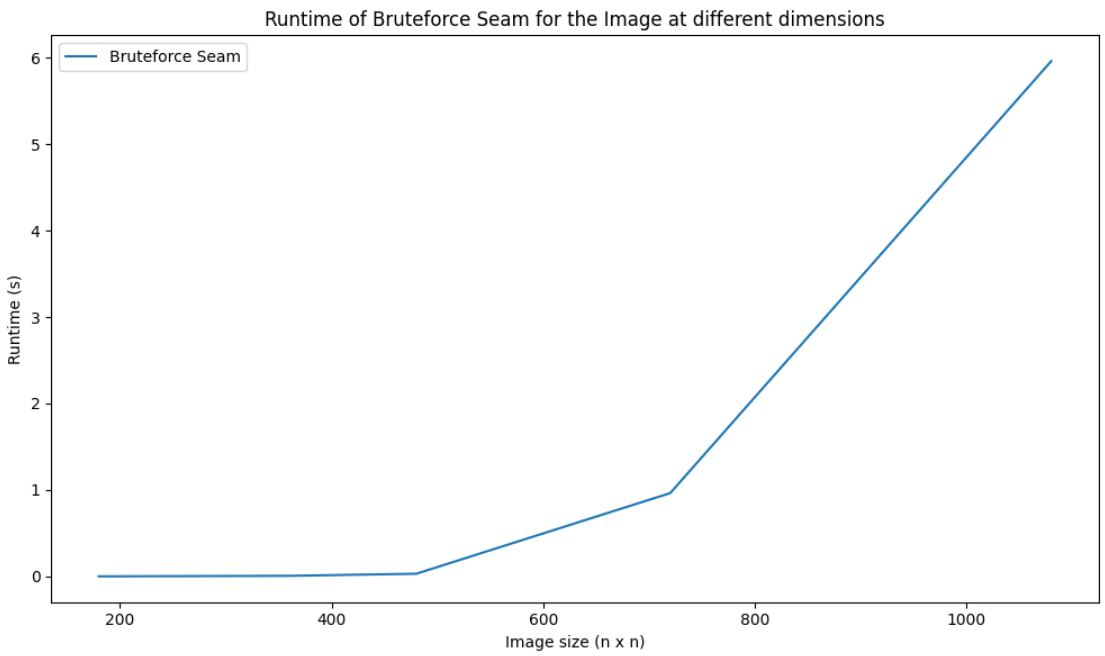}
    \caption{Time Complexity $T(n,n)$ of Minimum Seam using Naive algorithm}
    \label{fig:enter-label}
\end{figure}

\subsection{Greedy Algorithm}
Figure 9,10 shows the runtime for greedy algorithm for images of dimensions. The results of the plot are the same as our theoretical analyses i.e $C(n,n,s) = O(n^2)$ and $T(n,n) = O(n)$

\begin{figure}[hbtp]
    \centering
    \includegraphics[width=0.75\linewidth]{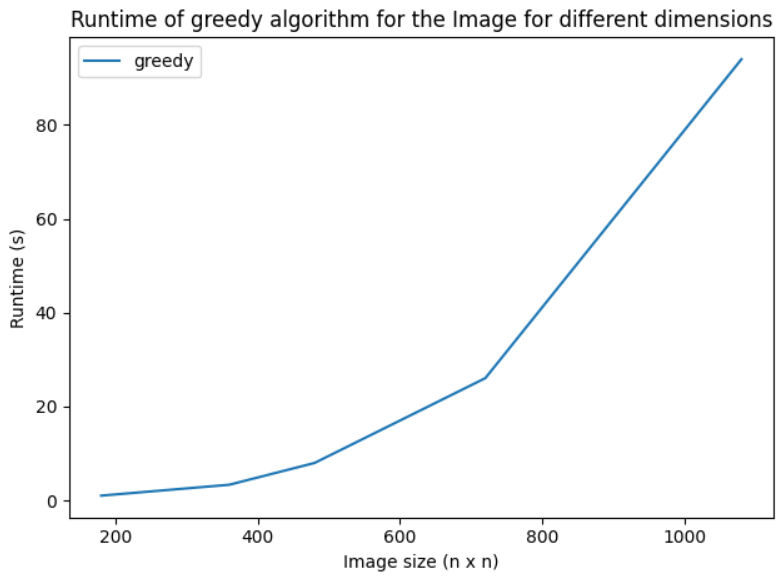}
    \caption{Overall time complexity of the greedy algorithm}
    \label{fig:enter-label}
\end{figure}

\begin{figure}[hbtp]
    \centering
    \includegraphics[width=0.75\linewidth]{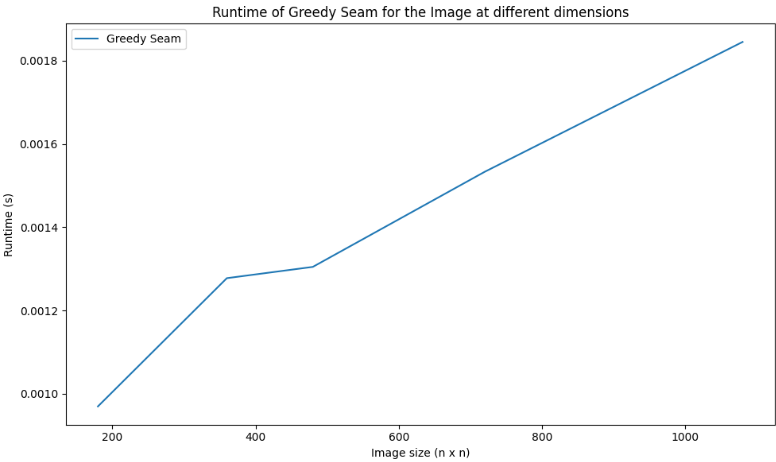}
    \caption{Time Complexity $T(n,n)$ of Minimum Seam using Greedy algorithm}
    \label{fig:enter-label}
\end{figure}


\subsection{Dynamic Algorithm}
Figure 11,12 shows the runtime for greedy algorithm for images of dimensions. The results of the plot are the same as our theoretical analyses i.e $C(n,n,s) = O(n^3)$ and $T(n,n) = O(n^2)$

\begin{figure}[hbtp]
    \centering
    \includegraphics[width=0.75\linewidth]{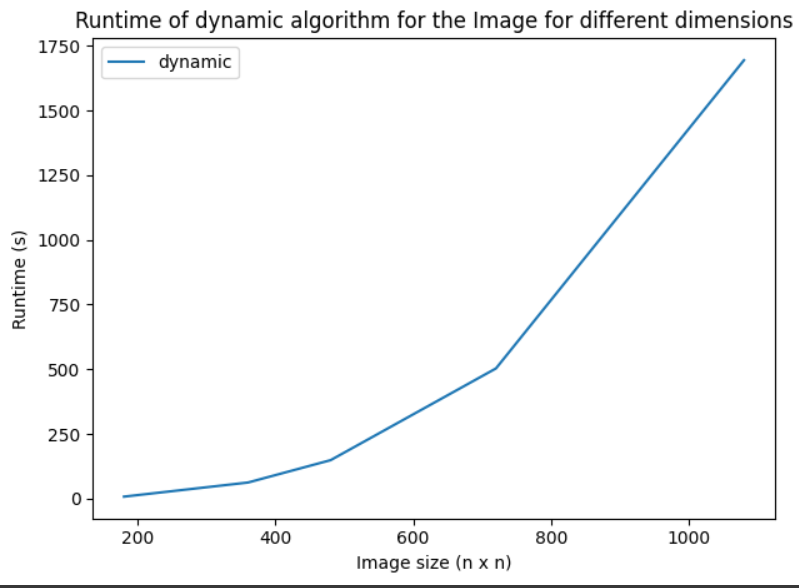}
    \caption{Overall time complexity of the dynamic algorithm}
    \label{fig:enter-label}
\end{figure}

\begin{figure}[hbtp]
    \centering
    \includegraphics[width=0.75\linewidth]{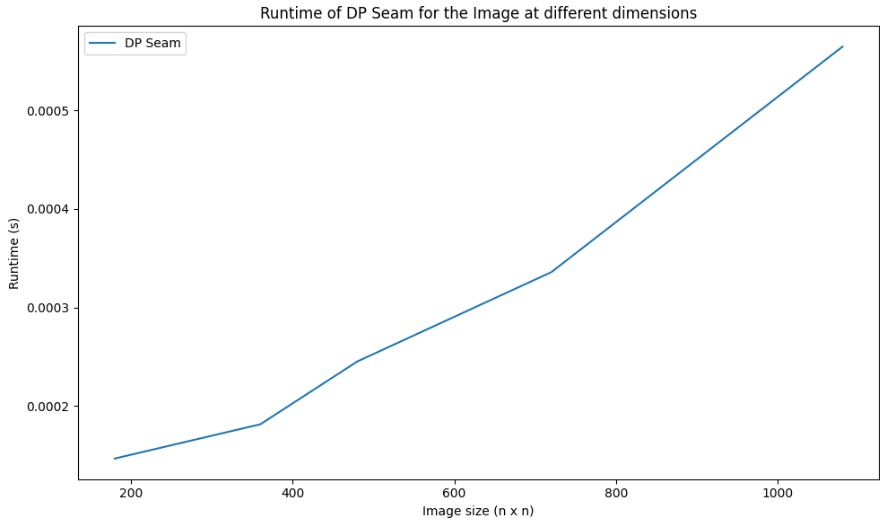}
    \caption{Time Complexity $T(n,n)$ of Minimum Seam using DP}
    \label{fig:enter-label}
\end{figure}
\subsection{GPU parallelized Algorithm}
Figure 13,14 shows the runtime for greedy algorithm for images of dimensions. The results of the plot are the same as our theoretical analyses i.e $C(n,n,s) = O(n^2)$ and $T(n,n) = O(n)$

\begin{figure}[hbtp]
    \centering
    \includegraphics[width=0.75\linewidth]{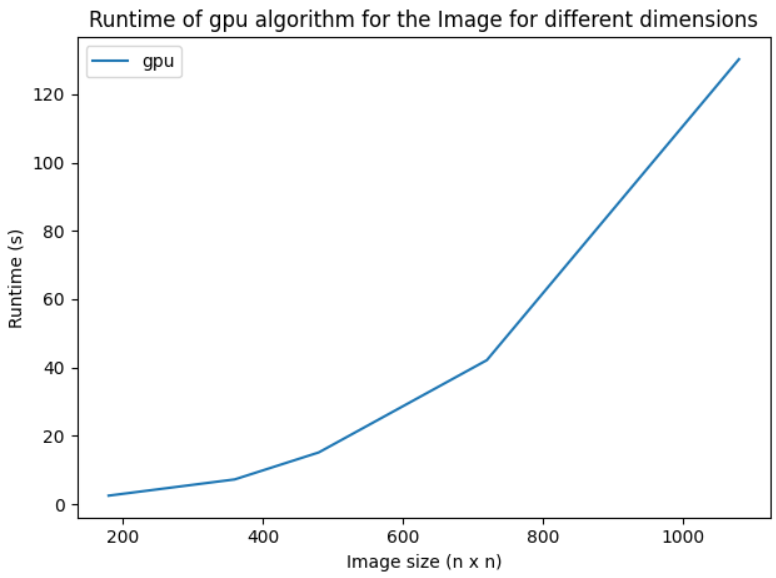}
    \caption{Overall time complexity of the GPU Dynamic algorithm}
    \label{fig:enter-label}
\end{figure}
\begin{figure}
    \centering
    \includegraphics[width=0.75\linewidth]{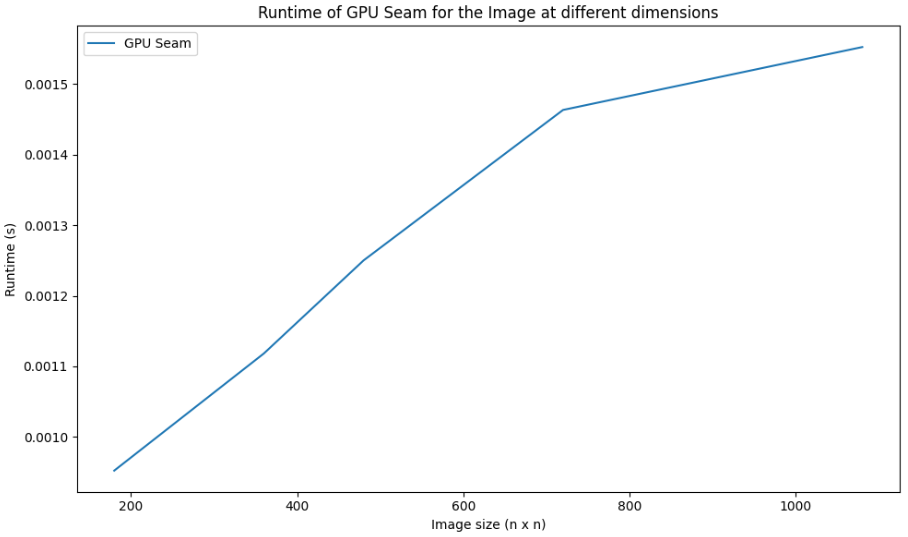}
    \caption{Time Complexity $T(n,n)$ of Minimum Seam using GPU}
    \label{fig:enter-label}
\end{figure}

\subsection{Overall Comparision}
The figures show an overall comparision of the greedy, dynamic and GPU based approach. We donot include the bruteforce algorithm here since it is not practical to run it on the image sizes in the plot.

\begin{figure}[hbtp]
    \centering
    \includegraphics[width=0.75\linewidth]{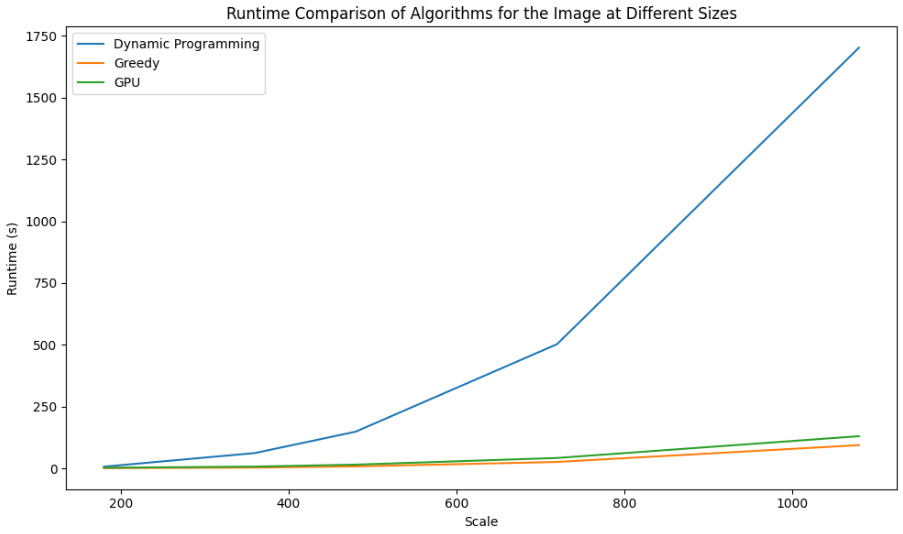}
    \caption{Overall Runtimes of dynamic, gpu and greedy solutions}
    \label{fig:enter-label}
\end{figure}

\begin{figure}][hbtp]
    \centering
    \includegraphics[width=0.75\linewidth]{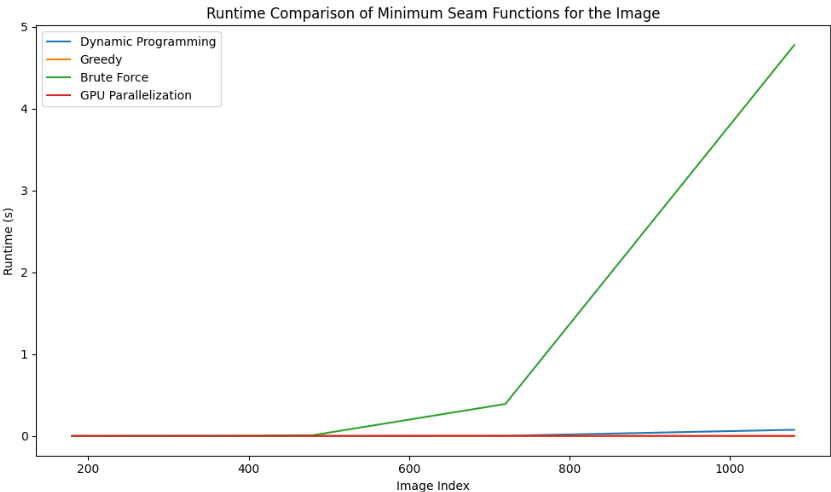}
    \caption{Runtime comparisions of $T(n,n)$ for different techniques}
    \label{fig:enter-label}
\end{figure}

\section{Results}
Here are some experimental results from the algorithm. These results are calculated using the GPU seam carving algorithm. These show an original image and its resized version. 

\begin{figure}%
    \centering
    \subfloat[\centering Original]{{\includegraphics[width=5cm]{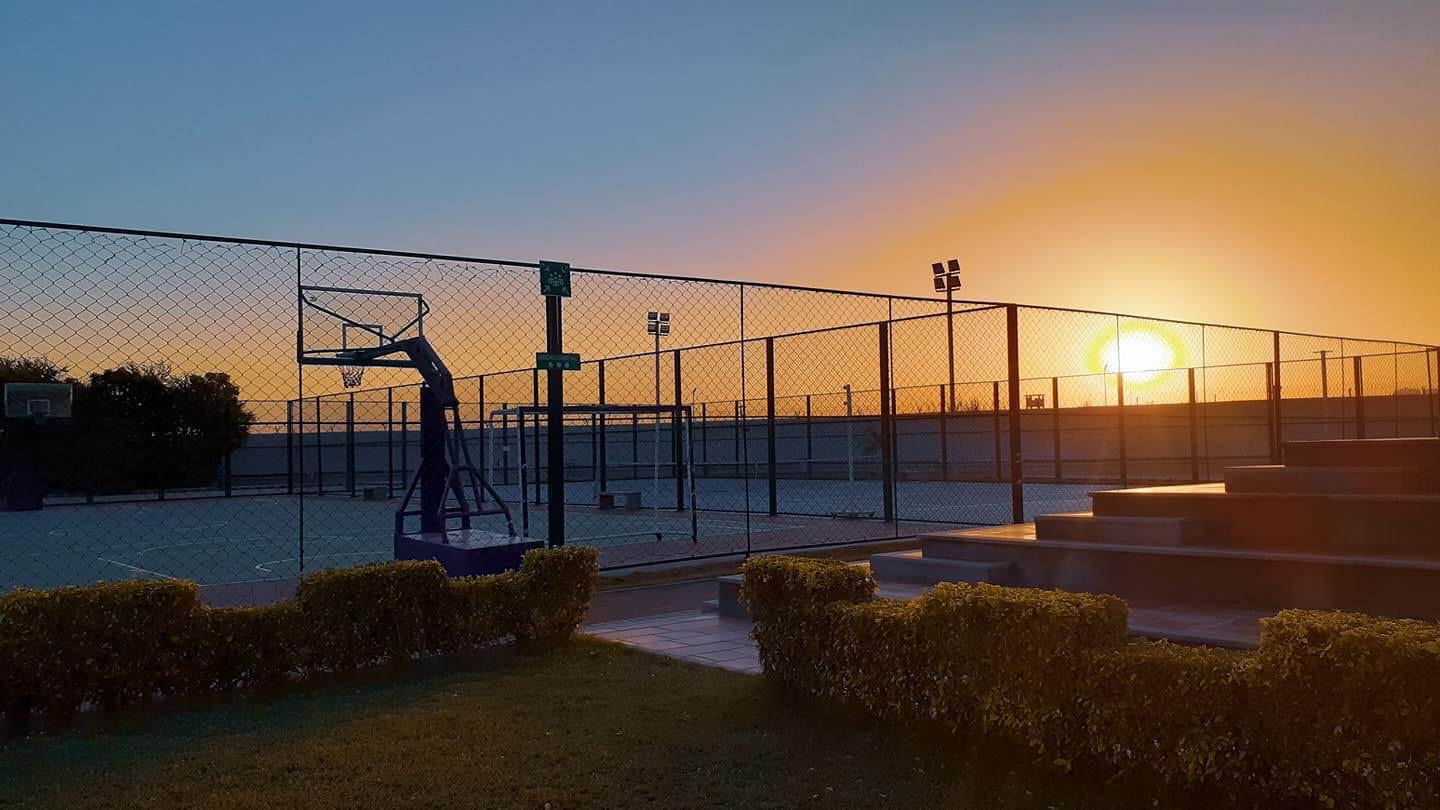} }}%
    \qquad
    \subfloat[\centering Scaled to 50\%]{{\includegraphics[width=5cm]{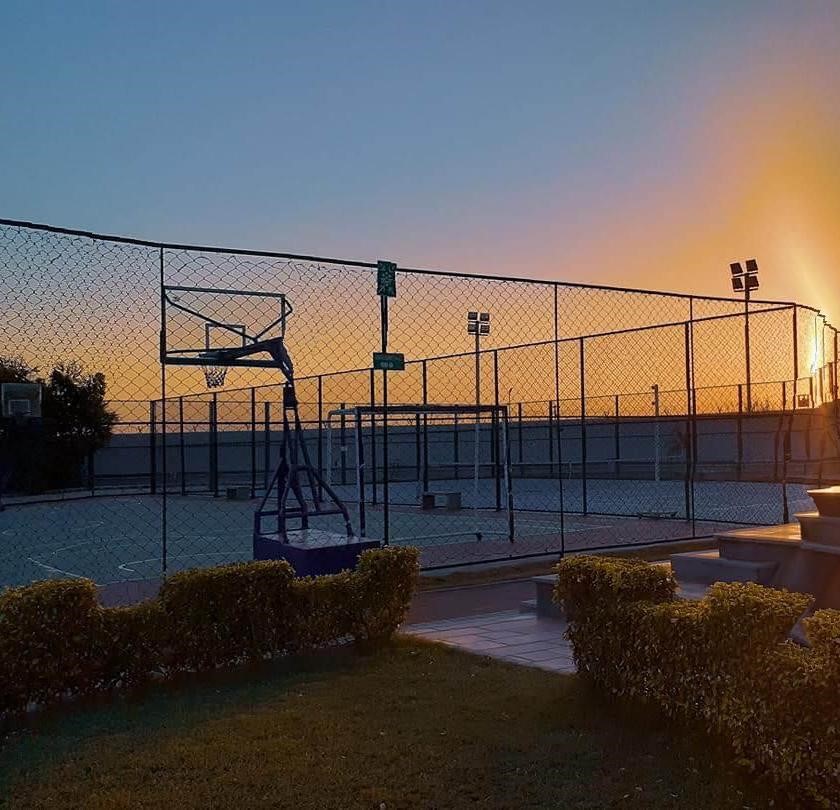} }}%
    \caption{Habib University courts landscape scaled using seam carving}%
    \label{fig:example}%
\end{figure}

\section{Conclusion}
Seam Carving is a technique used for image/video compression, enlargement, targeting and masking/object removal. This paper presents an analysis of different design techniques of the seam carving algorithm as well as validates their runtime complexities thorugh empirical experimentation.

\bibliographystyle{ACM-Reference-Format}
\bibliography{sample-sigplan}

\end{document}